\documentclass[fleqn]{mat01}
\usepackage{times,mathtimy,amssymb,latexsym}
\begin{document}

\setcounter{page}{71} \firstpage{71}

\renewcommand\theequation{\thesection\arabic{equation}}

\newtheorem{theore}{\bf Theorem}
\renewcommand\thetheore{\arabic{section}.\arabic{theore}}
\newtheorem{lem}[theore]{\it Lemma}
\newtheorem{coro}{\rm COROLLARY}
\newtheorem{propo}[theore]{\rm PROPOSITION}
\newtheorem{definit}[theore]{\rm DEFINITION}
\newtheorem{rem}[theore]{\it Remark}
\newtheorem{pot}[theore]{\it Proof of Theorem}

\newtheorem{exampl}[theore]{\it Example}

\title{Semigroups on Frechet spaces and equations with infinite delays}

\markboth{T Sengadir}{Semigroups on Frechet spaces}

\author{T SENGADIR}

\address{Department of Mathematics, SSN College of Engineering,
Old Mahabalipuram Road, Kalavakkam 603 110, India}

\volume{117}

\mon{February}

\parts{1}

\pubyear{2007}

\Date{MS received 28 February 2005}

\begin{abstract}
In this paper, we show existence and uniqueness of a solution to a
functional differential equation with infinite delay. We choose an
appropriate Frechet space so as to cover a large class of
functions to be used as initial functions to obtain existence and
uniqueness of solutions.
\end{abstract}


\maketitle

\section{Introduction and preliminaries}

In this paper we study linear functional differential equations
with infinite delay. Consider
\begin{align}
x'(t) &=ax(t) +\sum\limits_{i=1}^{\infty} b_{i}x(t
-\tau_{i}),\quad t
\geq 0\nonumber\\[.4pc]
x(\theta) &=\phi(\theta),\quad \theta  {\in}(-\infty,0]
\end{align}
where $a \in \mathbb{R}, \{b_{i}\}^{\infty}_{i=1}$ is an arbitrary
sequence of real numbers, $\{\tau_{i}\}^{\infty}_{i=1}$ is a
strictly increasing sequence of strictly positive reals such that
$\hbox{lim}_{i \rightarrow \infty} \tau_{i} = \infty$ and
$\phi\hbox{\rm :}\  (-\infty, 0]\longrightarrow\mathbb{R}$ is
continuous.

For the special case $\{b_{i}\}^{\infty}_{i=1} \in l^{1}$, (1.1)
can be uniquely solved for any given $\phi\in$ BC $(-\infty,0]$,
the space of all bounded real-valued continuous functions. The
proof of this is indicated in Example1.2. Denote this solution by
$x_{\phi}$. Consider the family of operators $S_{t}, t \geq 0$ on
BC $(-\infty,0]$ defined as
\begin{align}
[S_{t}\phi](\theta) &= x_{\phi} (t + \theta), \quad\hbox{if}\quad t + \theta >
0\nonumber\\[.2pc]
&= \phi(t + \theta), \quad\hbox{if}\quad t + \theta \leq 0.
\end{align}
It is elementary to see that $\{S_{t}\}$ is {\it not} a strongly
continuous semigroup on BC $(-\infty, 0]$. We prove this as
follows\hbox{\rm :}\  If $S_{t}$ is a semigroup, we must have
\begin{equation*}
\lim\limits_{t\rightarrow 0} S_{t}\phi =\phi.
\end{equation*}
Given $\epsilon > 0$, we can find $\delta > 0$ such that
$\|S_{t}\phi -\phi\|_\infty \leq \epsilon$ for $|t|\leq \delta$.
Now, let $\delta^{*} = \hbox{min} (1,\delta)$. Consider
$\theta_{1}, \theta_{2} \in (-\infty,-1]$ with $0 \leq \theta_{2}
-\theta_{1}\leq \delta^{*}$. Let $t = \theta_{2} -\theta_{1}$.
Note that $0 \leq t \leq \delta$ and
\begin{equation*}
\theta_{1} + t\leq \theta_{2} + \delta \leq \theta_{2} +1 \leq 0.
\end{equation*}
Consider
\begin{align*}
|\phi(\theta_{2}) -\phi(\theta_{1})| &\leq |\phi (t +\theta_{1})
-\phi(\theta_{1}|\\[.2pc]
&\leq |(S_{t}\phi) (\theta_{1}) - \phi(\theta_{1})|\\[.2pc]
&\leq \|S_{t}\phi -\phi\|_\infty\\[.2pc]
&\leq \epsilon.
\end{align*}
Thus, we have shown that $\phi$ is uniformly continuous on
$(-\infty, -1]$. By the uniform continuity of $\phi$ on $[-1, 0]$,
uniform continuity of $\phi$ on $(-\infty, 0]$ follows. But this
is a contradiction as there are bounded continuous functions which
are not uniformly continuous.

On the space BUC $(-\infty, 0]$, the space of bounded uniformly
continuous functions, the family of operators $S_{t}$ defined by
(1.2) do form a semigroup but this space is properly contained in
BC$(-\infty, 0]$. In the literature, certain Banach spaces, which
contain BC$(-\infty, 0]$, which are classes of functions
satisfying certain growth conditions, are used in the context of
infinite delay equations. Refer~\cite{4,5}. Our approach is to
find a Frechet space that contains BC$(-\infty, 0]$ on which
$S_{t}\hbox{'s}$ form a semigroup.

Instead of constructing a weight function which is related to
$b_{i}$ and $\tau_{i}$, we obtain a family of semi-norms for the
initial function $\phi$ that enable us to get estimates for the
solution and also to capture a Frechet space $F$ such that (1.2)
defines a strongly continuous semigroup on $F$. We do not make any
explicit summability assumption on the sequence
$\{b_{i}\}^{\infty}_{i=1}$, but the space $F$ heavily depends on
the properties of $\{b_{i}\}^{\infty}_{i=1}$. If $\{b_{i}\} \in
l^{1}$, then BC$(-\infty,0] \subset F$ but if $b_{i} =
\frac{1}{i}$, BC$(-\infty,0]$ is not contained in $F$.

The basic theory of finite delay differential equations is covered
in \cite{1}. \cite{2} and \cite{3} are some basic references for
Banach phase spaces related to infinite delay equations. Consider
the following examples\hbox{\rm :}\

\begin{exampl}
{\rm Let $X$ = BC$[0,\infty)$ and $A\hbox{\rm :}\  {\mathcal D}(A)
= \hbox{BC}^{1}[0,\infty)\longrightarrow \hbox{BC}[0,\infty)$ be
defined as $[A\varphi](x) ={\varphi}'(x)$. Define $S_{t}\hbox{\rm
:}\  X\longrightarrow X, t \geq 0$, as $(S_{t}\varphi)(x) =
\varphi(t + x)$. It is easy to see that $A$ satisfies the
following conditions\hbox{\rm :}\
\begin{enumerate}
\renewcommand\labelenumi{{\rm(\alph{enumi})}}
\leftskip .15pc

\item the resolvents $(\lambda{I}-A)^{-1}u =\int^{\infty}_{0}
\hbox{e}^{-\lambda{\tau}}S_{\tau}u\hbox{d}\tau$ exist for $\lambda
> 0,$ where $I$ is the identity operator on $X,$

\item $\|(\lambda{I} - A)^{-1}\|\leq\frac{1}{\lambda}$ and that
$S_{t}$ satisfy the following conditions\hbox{\rm :}\vspace{-.3pc}

\leftskip .4pc

\vspace{6pt} (i) $S_{0} = I,$

(ii) $S_{t + s} = S_{t}S_{s}$.

Despite the above observations $S_{t}$ is not strongly semigroup
on $X$ since the condition

(iii) $\lim_{t \rightarrow 0} S_{t}\varphi = \varphi, {\varphi}
\in X$

does {\it not} hold.\vspace{-1.2pc}

\end{enumerate}

The proof that (iii) does not hold is similar to the proof of the
analogous assertion in the case of the infinite delay equations.
The Hille--Yosida theorem is not applicable precisely because $A$
is not densely defined. One way to overcome this difficulty is to
construct a {\it smaller} Banach space called the Hille--Yosida
space for the operator $A$ \cite{11} on which the restriction of
$A$ generates a semigroup. In this construction, the condition (a)
plays a crucial role.

But consider the operator ${\mathcal A}\hbox{\rm :}\  {\mathcal
D}({\mathcal A}) = {\bf C}^{1}[0,\infty)\longrightarrow {\bf
C}[0,\infty)$ defined as $[{\mathcal A} \varphi](x)
={\varphi}'(x)$. ${\mathcal A}$ is an `extension' of $A$ that does
generate the semigroup ${\mathcal S}_{t}\hbox{\rm :}\ {\bf
C}[0,\infty) \longrightarrow {\bf C}[0,\infty)$ defined as
$[{\mathcal S}_{t}]\varphi(x)=\varphi(t + x)$ on the {\it Frechet
space} ${\bf C}[0,\infty)$. Thus$,$ it is clear that by
considering the Hille--Yosida space$,$ a lot of useful information
is lost.}
\end{exampl}

\begin{exampl}
{\rm We now indicate the proof of assertion that for $\{b_{i}\}
\in l^{1}$, (1.1) can be uniquely solved for any given $\phi \in
\hbox{BC}(-\infty, 0]$.

For $t \in[0,\tau_{1}], t -\tau_{i} \in (-\infty, 0]$ and hence
for $\phi \in\hbox{BC}(-\infty,0], \phi(t -\tau_{i})$ is
meaningful. Define $z_{i} \in C[0,\tau_{1}]$ as $z_{i} =
b_{i}\phi(t -\tau_{i})$. Now $\|z_{i}\|_\infty = |b_{i}|
\sup_{t\epsilon [0,\tau_{1}]} |\phi(t -\tau_{i})| \leq
|b_{i}|\|\phi\|_\infty = |b_{i}|\|\phi\|_\infty$. Since
$\{b_{i}\}$ is in $l^{1}$, the series $\sum_{i=1}^{\infty}
\|z_{i}\|$ converges. Hence $\sum_{i=1}^{\infty}z_{i}$ converges
in $C[0,\tau_{1}]$. Thus, $\sum_{i=1}^{\infty}\phi (t -\tau_{i})
\in C[0,\tau_{1}]$. Now, consider the ordinary differential
equation
\begin{align*}
x'(t) &= ax(t) + \sum\limits_{i=1}^{\infty} b_{i}\phi (t
-\tau_{i}), \quad t \in [0,\tau_{1}],\\[.4pc]
x(0) &=\phi(0)
\end{align*}
whose solution exists and is unique.

Now, assume the existence of a unique function $y_{k}\hbox{\rm :}\
(-\infty, k \tau_{1}] \longrightarrow\mathbb{R}$ such that
$y_{k}(\theta) = \phi(\theta)$ for $\theta {\in} (-\infty, 0]$ and
whose restriction to $[0, k\tau_{1}]$ is a solution to
\begin{align*}
x'(t) &= ax(t) + \sum\limits_{i=1}^{\infty} b_{i}x(t
-\tau_{i}),\quad t \epsilon [0, k\tau_{1}],\\[.4pc]
x(0) &=\phi(0).
\end{align*}

As before, the expression $\sum^{\infty}_{i=1} b_{i}y_{k} (t
-\tau_{i})$ defines a continuous function on $[k \tau_{1}, (k
+1)\tau_{1}]$. Now, consider
\begin{align*}
x'(t) &=ax(t) + \sum\limits_{i=1}^{\infty} b_{i}y_{k} (t
-\tau_{i}),\quad
t \in [k \tau_{1}, (k + 1)\tau_{1}],\\
x(k \tau_{1}) &=y_{k}(k\tau_{1})
\end{align*}
which has a unique solution $z$. Define $y_{k +1}$ as $y_{k +
1}(s)= y_{k}(s), s \in (-\infty, k\tau_{1}]$ and $y_{k+1}(s) =
z(s), s \in [k\tau_{1}, (k +1) \tau_{1}]$. Thus, we obtain a
solution to (1.1) in $(-\infty, (k +1)\tau_{1}]$. By induction and
patching up of solutions, we get a unique solution to (1.1) on the
whole of $(-\infty, \infty)$.

Let $X =BC(-\infty,0]$ and $\{b_{i}\}^{\infty}_{i=1} \in l^{1}$.
Define
\begin{align*}
A\hbox{\rm :}\  {\mathcal D}(A) =\left\{\varphi \in \hbox{BC}^{1}
(-\infty, 0] \hbox{\rm :}\ {\varphi}'(0)
=\sum\limits^{\infty}_{i=1}
b_{i}\varphi(-\tau_{i})\right\}\longrightarrow
{\hbox{BC}}(-\infty,0]
\end{align*}
as
\begin{equation*}
[A\varphi](x) = {\varphi}'(x).
\end{equation*}
${\mathcal D}(A)$ is not dense in $X$ and hence $A$ does not
generate a semigroup. Motivated by Example~1.1, we look for a
Frechet space $F$ that contains BC$(-\infty, 0]$ and an
`extension' ${\mathcal A}$ to $A$ such that ${\mathcal A}$
generates a semigroup on $F$.}
\end{exampl}

\begin{rem}
{\rm  The general theory of semigroups on Frechet spaces is very
complicated. For example$,$ even a bounded linear operator on a
Frechet space need not generate a semigroup \cite{13}. In
\cite{12}$,$ a generalisation of the Hille--Yosida theorem for a
closed and unbounded operator in a locally convex space is proved.
But the hypotheses of this theorem are not easily verified in many
concrete cases. By proving various estimates for the solution $x$
of (1.1)$,$ we are able to capture a Frechet space $F$ on which
the solution of (1.1) gives rise to a semigroup.

Refer to \cite{7}, \cite{14} and \cite{6} for applications of
semigroups on locally convex spaces to PDE's.

We need the following definitions and results in the next
section.}
\end{rem}

\begin{definit}$\left.\right.$\vspace{.5pc}

\noindent {\rm
\begin{enumerate}
\renewcommand\labelenumi{{\rm (\roman{enumi})}}
\leftskip .15pc
\item A topological vector space $X$ is said to be a Frechet space
if its topology is generated by a family of countable semi-norms
$\{q_{i}\}_{i=1}^{\infty}$ and $X$ is complete with respect to the
family $\{q_{i}\}^{\infty}_{i=1}$.

\item A linear map $S\hbox{\rm :}\  X \rightarrow X$ is said to be bounded if for
every $i \in {\mathbb N}$, there are finitely many indices $j_{1},
j_{2}, \dots, j_{m}$ and a constant $C$ such that for all $\phi
\in X$,
\begin{equation*}
\hskip -1.25pc q_{i}(S\phi) \leq C \max (q_{j_{1}}(\phi),
q_{j_{2}}(\phi),\dots ,q_{j_{m}}(\phi)).
\end{equation*}
\end{enumerate}
The basic theory of Frechet spaces and the proof of the following
proposition can be found in \cite{10}.}
\end{definit}

\begin{propo}$\left.\right.$\vspace{.5pc}

\noindent A linear map $S\hbox{\rm :}\  X\rightarrow X$ is
continuous if and only if it is bounded.
\end{propo}

\begin{definit}$\left.\right.$\vspace{.5pc}

\noindent {\rm A family of bounded linear operators
$\{S_{t}\hbox{\rm :}\ t \geq 0\}$ on $X$ is said to be a strongly
continuous semigroup if the properties (i), (ii) and (iii) of
Example~1.1 hold.}
\end{definit}

\begin{propo}$\left.\right.$\vspace{.5pc}

\noindent Let $l^{1}(X)$ denote the Banach space of all sequences
$\{x_{i}\}^{\infty}_{i =1}$ of elements of a Banach $X$ such that
$\sum^{\infty}_{i=1} \|x_{i}\| < \infty$. Let ${\bf C}([0, T];
l^{1})$ be the Banach space of all continuous functions
$h\hbox{\rm :}\  [0, T]\longrightarrow l^{1}$. The Banach space
$l^{1}({\bf C}([0, T]))$ is isometrically embedded in ${\bf C}([0,
T];l^{1})$.
\end{propo}

\begin{proof}
Let $G \in l^{1} ({\bf C}([0,T]))$. Define $G^{*}\hbox{\rm :}\
[0,T] \longrightarrow l^{1} (\mathbb{R})$ as $[G^{*}(s)]^{(i)} =
G^{(i)}(s)$. It is given that
\begin{equation*}
\sum\limits_{i =1}^{\infty} \hbox{sup} \{|G^{(i)} (s)|\hbox{\rm
:}\  s \in [0, T]\} < \infty.
\end{equation*}
\end{proof}

It is easy to check that $G^{*}$ is a bounded function and that
\begin{equation*}
\hbox{sup} \left\{\sum\limits_{i
=1}^{\infty}|[G^{*}(s)]^{(i)}|\hbox{\rm :}\  s \in [0,T]\right\} =
\sum\limits_{i =1}^{\infty} \hbox{sup} \{|G^{(i)}(s)|\hbox{\rm :}\
s \in [0,T]\}.
\end{equation*}
It remains to be shown that
\begin{equation*}
\mathop{\hbox{lim}}\limits_{\substack{t \rightarrow
t_{0}}}\sum\limits_{i =1}^{\infty}|[G^{*}(t)]^{(i)} -
[G^{*}(t_{0})]^{(i)}| = 0.
\end{equation*}
By the hypothesis for each $i$, $\hbox{lim}_{t\rightarrow
t_{0}}|[G^{*}(t)]^{(i)} - G^{*}(t_{0})]^{(i)}| = 0$ and for a
given $\epsilon > 0$, there exists $K \in \mathbb{N}$ such that
\begin{equation*}
\sum\limits^{\infty}_{i =K + 1} \hbox{sup} \{|G^{(i)}
(s)|\hbox{\rm :}\  s \in [0, T]\} < \epsilon/3.
\end{equation*}

Further, there exists $\delta > 0$ such that $|t -t_{0}| < \delta$
implies that
\begin{equation*}
\sum\limits^{K}_{i =1}|[G^{*}(t)]^{(i)} - [G^{*}(t_{0})]^{(i)}| <
\epsilon /3.
\end{equation*}
Hence for all $t, t_{0}$ in $[0, T]$ with $|t - t_{0}| < \delta$,
\begin{align*}
\sum\limits^{\infty}_{i =1}|[G^{*}(t)]^{[(i)} \!-\!
[G^{*}(t_{0})]^{[(i)}| &= \sum\limits^{K}_{i =1}|[G^{*}(t)]^{[(i)}
- [G^{*}(t_{0})]^{[(i)}|\\[.5pc]
&\quad\,+ \sum\limits^{\infty}_{i =K +1}|[G^{*}(t)]^{[(i)} -
[G^{*}(t_{0})]^{[(i)}| \leq {\epsilon}/3\\[.5pc]
&\quad\,+ \sum\limits^{\infty}_{i =K + 1}|[G^{*}(t)]^{[(i)}| +
\sum\limits^{\infty}_{i =K + 1}|[G^{*}(t_{0})]^{[(i)}|\\[.5pc]
&\leq {\epsilon}/3 + \epsilon/3 + \epsilon/3 = \epsilon.
\end{align*}
The result is proved.

The following definitions and statements on Frechet space valued
Riemann integral are found in \cite{9}.

\begin{theore}[\!]
Let $X$ be a Frechet space and let $u\hbox{\rm :}\  [a,b]
\longrightarrow X$ be continuous. The integral $\int^{b}_{a}u(t)
\hbox{\rm d}t \in X$ can be defined uniquely which has the
following properties\hbox{\rm :}\
\begin{enumerate}
\renewcommand\labelenumi{{\rm (\roman{enumi})}}
\leftskip .4pc
\item for every continuous linear functional $x^{*}\hbox{\rm :}\  X
\longrightarrow\mathbb{R}${\rm ,}
\begin{equation*}
\hskip -1.2pc x^{*} \left(\int^{b}_{a} u(t)\hbox{\rm d}t\right)
=\int^{b}_{a}x^{*}(u(t))\hbox{\rm d}t,
\end{equation*}

\item for every seminorm $q_{k}${\rm ,}
\begin{equation*}
\hskip -1.2pc q_{k} \left(\int^{b}_{a}u(t)\hbox{\rm d}t\right)\leq
{\int}^{b}_{a} q_{k}(u(t))\hbox{\rm d}t\hbox{\rm ,}
\end{equation*}
\item $\int^{b}_{a}u(t)\hbox{\rm d}t + \int^{c}_{b}u(t)\hbox{\rm d}t
=\int^{c}_{a}u(t)\hbox{\rm d}t$\hbox{\rm ,}\pagebreak
\item $\int^{b}_{a}[u(t) + v(t)]\hbox{\rm d}t = \int^{b}_{a}u(t)\hbox{\rm d}t
+\int^{b}_{a}v(t)\hbox{\rm d}t$\hbox{\rm ,}
\item $c \int^{b}_{a}u(t)\hbox{\rm d}t = \int^{b}_{a}cu(t)\hbox{\rm d}t$.
\end{enumerate}
\end{theore}

\begin{definit}$\left.\right.$\vspace{.5pc}

\noindent {\rm Let $J $ be a sub interval of $\mathbb{R}$. A
function $u\hbox{\rm :}\  J \longrightarrow X$ is said to be
differentiable at $t_{0} \in J$ if there exists $y \in X $ with
the following property\hbox{\rm :}\  for every $\epsilon > 0$ and
$k\in \mathbb{N}$, there exists $\delta > 0$ such that for all $h
\in \mathbb{R}$ with $|t_{0} + h| < \delta$ and $t_{0} + h \in J$,
we have
\begin{equation*}
q_{k} \left(\frac{u(t_{0} + h) - u(t_{0})} {h} - y\right) <
\epsilon.
\end{equation*}
Further, the element $y$ is denoted by $u'(t_{0})$.}
\end{definit}

\begin{definit}$\left.\right.$\vspace{.5pc}

\noindent {\rm Let $J$ be a sub interval of $\mathbb{R}$.
A~function $u\hbox{\rm :}\  J\longrightarrow X$ is said to be
continuously differentiable on $J$, if $u$ is differentiable at
every point $t_{0} \in J$ and the function mapping $t$ to $u'(t)$
is continuous on $J$. The class of all such functions is denoted
by ${\bf C}^{1}(J;X)$.}
\end{definit}

\section*{\it Fundamental theorem of integral calculus}

If $u\hbox{\rm :}\  [a,b]\longrightarrow X$ is continuously
differentiable, then
\begin{equation*}
u(a) - u(b) =\int^{b}_{a}u'(t){\rm d}t.
\end{equation*}

\section*{\it A Frechet phase space}

For a given $\varphi \in {\bf C}((-\infty, 0])$ define the family
of seminorms $\{\|.\|_{k}\hbox{\rm :}\  k \in \mathbb{N}\}$ as
\begin{equation*}
\|\varphi\|_{k} = \sup \{|\varphi(\theta)|\hbox{\rm :}\  \theta
\in [-k, 0]\}.
\end{equation*}
Let $\{b_{i}\}^{\infty}_{i =1}$ and $\{\tau_{i}\}^{\infty}_{i =1}$
be as in the beginning of this section. Given $k \in \mathbb{N}$,
define $n(k)\in \mathbb{N}$ as the smallest positive integer such
that $\tau_{i} \geq k\tau_{1}$ for all $i \geq n(k)$.

Define $F$ as
\begin{equation*}
F = \{\varphi \in {\bf C}(-\infty, 0]\hbox{\rm :}\ p_{k}(\varphi)
< \infty \quad\hbox{for all} \ \ k \in \mathbb{N}\},
\end{equation*}
where the seminorms $p_{k}$ are defined as follows\hbox{\rm :}\
\begin{equation*}
p_{k}(\varphi) = \sum\limits^{\infty}_{i = n(k)} \sup
\{|b_{i}\varphi(s -\tau_{i})|\hbox{\rm :}\ s \in [0,k\tau_{1}]\}.
\end{equation*}

\begin{propo}$\left.\right.$\vspace{.5pc}

\noindent The space $F$ equipped with the topology generated by
the family of seminorms $\{\|\cdot\|_{k}\hbox{\rm :}\ k \in
\mathbb{N}\} \cup \{p_{k}\hbox{\rm :}\  k \in \mathbb{N}\}$ is a
Frechet space.
\end{propo}

\begin{proof}
Let $\phi_{j}$ be a Cauchy sequence in $F$. Clearly, there exists
$\phi \in {\bf C}(-\infty, 0]$ such that $\phi_{j}$ converges to
$\phi$ uniformly on every compact set of $(-\infty, 0]$. Consider
the Banach space $l^{1}({\bf C}([0, k\tau_{1}]))$. For every $j
\in \mathbb{N}$, define $G_{j} \in l^{1} (C([0{\rm ,}
k\tau_{1}]))$ as
\begin{equation*}
G_{j}^{(i)}(s) = b_{n(k) +(i -1)}\phi_{j}(s -\tau_{n(k) +
(i-1)}){\rm ,}\quad i \in \mathbb{N}.
\end{equation*}
By the hypothesis, $G_{j}$ is a Cauchy sequence in $l^{1}({\bf
C}([0,k{\tau}_{1}]))$ and hence converges to some $G \in
l^{1}({\bf C}([0, k\tau_{1}]))$. For fixed $i$ and $s$,
$\lim_{j\rightarrow\infty} G^{(i)}_{j}(s) = G^{(i)}(s)$ and hence
$G^{(i)}(s) = b_{n(k) + (i -1)}\phi (s -\tau_{n(k) + (i -1)})$.
This implies that $p_{k}(\phi) < \infty$ and $\phi_{j}$ converges
to $\phi$ in $F$. The proof is complete.
\end{proof}

Now we state the definition of a mild solution of an abstract
Cauchy problem in Frechet spaces from \cite{8}.

Let $X$ be a Frechet space whose topology is given by the family
of seminorms $\{q_{k}\hbox{\rm :}\  k \in \mathbb{N}\}$ and
${\mathcal D}(A)$ be a subset of $X$. For a closed operator
$A\hbox{\rm :}\ {\mathcal D}(A) \longrightarrow X$ and $\phi \in
X$, consider the abstract Cauchy problem
\begin{align}
\frac{\hbox{d}u}{\hbox{d}t}& = Au,\nonumber\\[.4pc]
u(0) &=\phi.
\end{align}

$\left.\right.$\vspace{-1pc}

\begin{definit}$\left.\right.$\vspace{.5pc}

\noindent {\rm A function $u(\cdot) \in {\bf C}([0,\infty),
{\mathcal D}(A))\cap {\bf C}^{1}([0,\infty),X)$ that satisfies
(1.3) is said to be a {\it solution} of Problem (1.3). A function
$u(\cdot) \in {\bf C}([0,\infty),X)$ is said to be a {\it mild
solution} of (1.3), if $v(t) = \int^{t}_{0}u(s)\hbox{d}s \in
{\mathcal D}(A)$, for all $t \geq 0$, and
\begin{equation}
\frac{\hbox{d}(v(\cdot))}{\hbox{d}t}(t) =A(v(t)) + \phi, \quad t
\geq 0.
\end{equation}}
\end{definit}

\section{Main results}

Let $F$ be as in the previous section. Define ${\mathcal D}(A) =
\{\varphi \in F\hbox{\rm :}\ {\varphi}' \in F$ and ${\varphi}'(0)
= L\varphi\}$ and $A\hbox{\rm :}\  {\mathcal D}(A)\longrightarrow
F$ as $A\varphi ={\varphi}'$ where
\begin{equation*}
L{\varphi} = a{\varphi}(0) + \sum\limits^{\infty}_{i=1}
b_{i}\varphi(-\tau_{i}).
\end{equation*}\vspace{.1pc}

\setcounter{theore}{0}
\begin{definit}$\left.\right.$\vspace{.5pc}

\noindent {\rm We say that a function $x\hbox{\rm :}\ \mathbb{R}
\longrightarrow\mathbb{R}$ is a solution to (1.1) if the following
hold\hbox{\rm :}\

\begin{enumerate}
\renewcommand\labelenumi{{\rm (\roman{enumi})}}
\leftskip .4pc
\item $x$ is continuous and $x(\theta) = \phi(\theta)$ for all
$\theta {\in}(-\infty,0]$.

\item The restriction of $x$ to $[0,\infty)$ is continuously
differentiable.

\item $x'(t) = ax(t) +\sum^{\infty}_{i =1} b_{i}x(t -\tau_{i})$ for
all $t \geq 0$.\vspace{-.5pc}
\end{enumerate}}
\end{definit}

\begin{rem}
{\rm  Note that our definition of the solution does not imply that
$x$ is differentiable from the left at $t=0$.}\vspace{1pc}
\end{rem}

\begin{theore}[\!]
Let $F,$ $L$ and $A$ be as in the beginning of this section. Then
$A$ generates a strongly continuous semigroup $\{S_{t}\hbox{\rm
:}\ t\geq 0\}$ of bounded linear operators on $F$. Further$,$ for
a given $\phi \in F,$ the map $x\hbox{\rm :}\
\mathbb{R}\rightarrow \mathbb{R}$ defined as
\begin{align*}
x(t) &= \phi(t),\quad t\in(-\infty, 0],\\[.2pc]
x(t) &= [S_{t}\phi](0),\quad t\in (0,\infty)
\end{align*}
is a unique solution to $(1.1)$.
\end{theore}

Besides, fixing $\phi \in F$ and defining $u\hbox{\rm :}\ [0,
\infty)\rightarrow F$ as $u(t)=S_{t}\phi, u(\cdot)$ is a mild
solution to the abstract Cauchy problem (1.3).

We need the following lemmas to prove Theorem 2.3 and we actually
define the semi-group via the solution to (1.1).

\begin{lem}
Let $\phi \in F$. The problem $(1.1)$ has a unique solution
$x\hbox{\rm :}\  \mathbb{R}\rightarrow \mathbb{R}$.
\end{lem}

Further$,$ for each $k \in \mathbb{N}$ there exists a finite
subset $\Lambda$ of $\{\|.\|_{k}\hbox{\rm :}\  k\in
\mathbb{N}\}\cup \{p_{k}\hbox{\rm :}\  k\in \mathbb{N}\}$ and a
constant $C_{k}\geq 0$ such that \setcounter{equation}{0}
\begin{equation}
\sup \{|x(s)|\hbox{\rm :}\  s\in [0, k\tau_{1}]\}\leq  C_{k} \max
\{q (\phi)\hbox{\rm :}\  q\in \Lambda\}.
\end{equation}

\begin{proof}
Consider $t \in [0, \tau_{1}]$. Clearly, $t-\tau_{i}\leq
t-\tau_{1}$ for all $i \in \mathbb{N}$. Thus for $t \in [0,
\tau_{1}], t-t_{i}\leq 0$. Define $y_{1}\hbox{\rm :}\
[0,\tau_{1}]\rightarrow \mathbb{R}$ as the unique solution of the
initial value problem.
\begin{align}
x'(t)&=ax (t) + \sum\limits_{i=1}^{\infty} b_{i}\phi
(t-\tau_{i}),\nonumber\\[.4pc]
x(0) &= \phi (0).
\end{align}
Note that as $\phi \in F, t \rightarrow \sum_{i=1}^{\infty}
b_{i}\phi (t-\tau_{i})$ defines a continuous function on $[0,
\tau_{1}]$. We have
\begin{equation*}
y_{1}(t)=\phi (0) \hbox{e}^{at} + \hbox{e}^{at} \int_{0}^{t}
\hbox{e}^{-as} \left(\sum \limits_{i=1}^{\infty} b_{i}\phi
(s-\tau_{i})\right)\hbox{d}s.
\end{equation*}
Define $x_{1}\hbox{\rm :}\  (-\infty, \tau_{1}]\rightarrow
\mathbb{R}$ as
\begin{align*}
x_{1}(s) &= \phi (s),\quad s \in (-\infty, 0]\\[.2pc]
&= y_{1} (s),\quad s \in [0, \tau_{1}].
\end{align*}
In the remaining  part of the proof we shall assume that $a \neq
0$. The estimates for $a=0$ are easier to obtain. Clearly,
\begin{equation}
\sup \{| x_{1}(t)|\hbox{\rm :}\  t \in [0, \tau_{1}]\}\leq
\left(\sup\limits_{s\in [0,\tau_{1}]} \hbox{e}^{as}\right) \|\phi
\|_{1} + \left(\sup\limits_{s\in [0,\tau_{1}]}
\frac{\hbox{e}^{as}-1}{a}\right)p_{1} (\phi).
\end{equation} 
Here, note that for $r>0, \frac{{\rm e}^{ar}-1}{a}>0$ for all
$a\neq 0$.  Now we claim that for each $k \in \mathbb{N}$, there
exists a function $x_{k}\hbox{\rm :}\  (-\infty, k
\tau_{1}]\rightarrow \mathbb{R}$ with the following
properties\hbox{\rm :}\

\begin{enumerate}
\renewcommand\labelenumi{{\rm (\roman{enumi})}}
\leftskip .4pc
\item For each $t \in [0, k \tau_{1}]$,
\begin{equation*}
\hskip -1.25pc \sum_{i=1}^{\infty} b_{i} x_{k} (t-\tau_{i})
\end{equation*}
converges and this summation defines a continuous function on $[0,
k \tau_{1}]$.

\item $x_{k}$ is the unique solution to
\begin{align}
\hskip -1.25pc x' (t) &= ax (t) + \sum_{i=1}^{\infty} b_{i}
x(t-\tau_{i}),
\quad\hbox{for} \ \ t \in [0, k \tau_{1}]\nonumber\\
\hskip -1.25pc x(\theta) &= \phi (\theta),\quad\hbox{for} \ \ t
\in (-\infty,
0].
\end{align}

\item There exists a finite subset $\Lambda$ of
$\{\|\cdot\|_{k}\hbox{\rm :}\  k\in \mathbb{N}\}\cup
\{p_{k}\hbox{\rm :}\ k\in \mathbb{N}\}$ and a constant $C_{k}\geq
0$ such that
\begin{equation*}
\hskip -1.25pc \sup \{|x_{k}(s)|\hbox{\rm :}\  s\in [0,
k\tau_{1}]\}\leq  C_{k} \max \{q (\phi)\hbox{\rm :}\  q\in
\Lambda\}.
\end{equation*}
\end{enumerate}
We prove this by induction on $k$.  The case $k=1$ is already
proved.

Assuming that our claim is true for arbitrary $k \in \mathbb{N}$,
we show that the claim is true for $k + 1$. Define
$I_{k}=[(k-1)\tau_{1},k\tau_{1}]$ and $\|x_{k}\|_{I_{k}}$ as
\begin{equation*}
\|x_{k}\|_{I_{k}}=\sup \{|x_{k}(s)|\hbox{\rm :}\  s \in
[(k-1)\tau_{1}, k\tau_{1}]\}.
\end{equation*}

For $s \in I_{k+1}$ and $i \in \mathbb{N}, s-\tau_{i} \in
(-\infty, k\tau_{1}]$ and so we can consider the summation
\begin{equation*}
\sum_{i=1}^{\infty} b_{i} x_{k}(s-\tau_{i}) =
\sum_{i=1}^{n(k+1)-1} b_{i} x_{k}(s-\tau_{i})+
\sum_{i=n(k+1)}^{\infty} b_{i} x_{k}(s-\tau_{i}).
\end{equation*}
The first summation involves only finitely many terms. For $s \in
[k\tau_{1}, (k + 1) \tau_{1}]$ and $i\geq n (k+1), s-\tau_{i}<
s-(k + 1) \tau_{1}\leq 0$ and hence
\begin{equation*}
\sum_{i=n(k+1)}^{\infty} b_{i} x_{k}(s-\tau_{i})
=\sum_{i=n(k+1)}^{\infty} b_{i} \phi (s-\tau_{i}).
\end{equation*}
So, for $s \in I_{k + 1}$, the summation $\sum_{i=1}^{\infty}
b_{i}x_{k} (s-\tau_{i})$ defines a continuous function on
$I_{k+1}$.

For a given $k \in \mathbb{N}-\{1\}$, define $m(k)$ as the
smallest positive integer such that
\begin{equation*}
-m(k) < \min \{(k-1)\tau_{1}-\tau_{i}\hbox{\rm :}\  i = 1,
2,\dots, n(k)-1\}.
\end{equation*}
The following estimates follow from the definition of the
seminorms and the integers $m(k)$\hbox{\rm :}\
\begin{equation}
\sup\left\{\left|\sum_{i=n(k + 1)}^{\infty}
b_{i}x_{k}(s-\tau_{i})\right|\hbox{\rm :}\  s\in [k\tau_{1}, (k +
1)\tau_{1}]\right\} \leq p_{k + 1}(\phi).
\end{equation}
Moreover,
\begin{align}
&\sup\left\{\left|\sum_{i=n}^{n(k+1)-1}
b_{i}x_{k}(s-\tau_{i})\right|\hbox{\rm :}\  s\in [k\tau_{1},
(k+1)\tau_{1}]\right\}\nonumber\\[.5pc]
&\quad\,\leq \sum_{i=1}^{n(k+1)-1} |b_{i}|\times \max
\{\|\phi\|_{m(k+1)}, \|x_{k}\|_{I_{k}}\}.
\end{align}

Define $y_{k+1}\hbox{\rm :}\  [k\tau_{1},(k+1)\tau_{1}]\rightarrow
\mathbb{R}$ as
\begin{equation*}
y_{k+1}(t)=x_{k}(k\tau_{1})
\hbox{e}^{a(t-k\tau_{1})}+\hbox{e}^{at} \int_{k\tau_{1}}^{t}
\hbox{e}^{-as} \left(\sum_{i=1}^{\infty}b_{i}x_{k}
(s-\tau_{i})\right)\hbox{d}s.
\end{equation*}

From (2.5) and (2.6), we get the estimate
\begin{align}
&\sup\left\{\left|\sum_{i=1}^{\infty}
b_{i}x_{k}(s-\tau_{i})\right|\hbox{\rm :}\  s\in [k\tau_{1},
(k+1)\tau_{1}]\right\}\nonumber\\[.4pc]
&\quad\,\leq \left[\sum_{i=1}^{n(k+1)-1} |b_{i}|\times \max
\{\|\phi\|_{m(k+1)},
\|x_{k}\|_{I_{k}}\}+p_{k+1}(\phi)\right].
\end{align}

Defining $x_{k+1}\hbox{\rm :}\  (-\infty,
(k+1)\tau_{1}]\rightarrow \mathbb{R}$ as
\begin{align*}
x_{k+1}(s) &= x_{k}(s), s \in (-\infty, k\tau_{1}]\\[.2pc]
&= y_{k+1}(s), s \in (k\tau_{1},(k+1)\tau_{1}],
\end{align*}
our claims (i) and (ii) are proved. Now, we proceed to prove (iii)
for $k+1$.  From the definition of $y_{k+1}$ and the estimate
(2.7), we get, for $k \in \mathbb{N}$,
\begin{align}
&\sup\{|x_{k+1}(s)|\hbox{\rm :}\  s \in [k\tau_{1}, (k+1)
\tau_{1}]\}=\|x_{k+1}\|_{I_{k+1}}\nonumber\\[.4pc]
&\quad\,= \sup\{|y_{k+1}(s)|\hbox{\rm :}\  s \in [k\tau_{1}, (k+1)
\tau_{1}]\}\nonumber\\[.4pc]
&\quad\,\leq \left(\sup\limits_{s \in[0,k\tau_{1}]}
\hbox{e}^{as}\right)\|x_{k}\|_{I_{k}} + \left(\sup\limits_{s
\in[0,\tau_{1}]}\frac{\hbox{e}^{a(k+1)s} - 1}{{a}}\right)\nonumber\\[.4pc]
&\qquad\ \times \left[\sum_{i=1}^{n(k+1)-1} |b_{i}| \times \max
\{\|\phi\|_{m(k+1)},
\|x_{k}\|_{I_{k}}\}+p_{k+1}(\phi)\right].
\end{align}
Since
\begin{equation*}
\|x_{k}\|_{I_{k}}\leq \sup \{|x_{k}(s)|\hbox{\rm :}\  s \in [0,
k\tau_{1}]\}
\end{equation*}
the assertion (iii) for $k+1$ follows from the estimate (2.8) and
hence the assertion (iii) for $k$.

The solution to (1.1) is obtained by patching the functions
$x_{k}$. Uniqueness of $x$ now follows.
\end{proof}

\begin{lem}
Let $\varphi \in F$ and $x\hbox{\rm :}\  \mathbb{R}\rightarrow
\mathbb{R}$ be a continuous function such that $x(\theta)=\varphi
(\theta)$ for all $\theta \in (-\infty, 0]$. Define $u\hbox{\rm
:}\ [0, \infty)\rightarrow {\bf C}((-\infty,0])$ as
$[u(t)](\theta)=x(t+\theta)$. Then $u \in {\bf C} ([0,
\infty);F)$.
\end{lem}

\begin{proof}
Fix $t \in [0,j\tau_{1}]$. It is trivial to check that $u(t) \in
{\bf C} ((-\infty,0])$. We also have the estimate
\begin{align}
\|u(t)\|_{k} &=\sup \{\theta \in [-k,
0]\hbox{\rm :}\  |[u(t)](\theta)|\}\nonumber\\[.2pc]
&\leq \max (\|\varphi \|_{k}, \sup \{|x(t)|\hbox{\rm :}\  t \in
[0, j\tau_{1}\}).
\end{align}
\end{proof}

Now we show that $p_{k}(u(t))< \infty$ for each $k \in
\mathbb{N}$. Let $s\in [0, k\tau_{1}]$ and $t \in [0, j\tau_{1}]$.
It is clear that $n(k+j)\geq n(k)$ and we have the following
assertions:
\begin{equation*}
i\geq n(k)\Rightarrow (s-\tau_{i})\leq 0
\end{equation*}
and
\begin{equation*}
n(k)\leq i< n(k+j)\Rightarrow -(k+j) \tau_{1} \leq
t+s-\tau_{i}\leq j\tau_{1}.
\end{equation*}
Thus, for every $K\geq n(k+j)$, we have
\begin{align*}
\sum_{i=n(k)}^{K} b_{i}[u(t)](s-\tau_{i})
&=\sum_{i=n(k)}^{n(k+j)-1} b_{i}[u(t)](s-\tau_{i})\\[.5pc]
&\quad\,+ \sum_{i=n(k+j)}^{K} b_{i}[u(t)](s-\tau_{i})\\[.5pc]
&=\sum_{i=n(k)}^{n(k+j)-1} b_{i}x(t+s-\tau_{i})\\&\quad\,+
\sum_{i=n(k+j)}^{K} b_{i}\phi(t+s-\tau_{i}).
\end{align*}
Now choosing a positive integer $m\geq (k+j)\tau_{1}$, we have the
estimate
\begin{equation}
p_{k}(u(t))\leq \sum_{i=1}^{n(k + j)}|b_{i}|\max
(\sup\{|x(t)|\hbox{\rm :}\  t \in [0, j\tau_{i}\},\|\phi\|_{m}) +
p_{(k + j)}(\phi).
\end{equation}
Equations (2.9) and (2.10) show that $u(t)\in F$.

Next, we show the continuity of $u$ at an arbitrary $t_{0} \in
[0,\infty)$. Let $\epsilon > 0$  be given. By the uniform
continuity of $x$ on $[-k-1 + t_{0},t_{0} + 1]$, there exists
$\delta > 0$ such that for $p,q \in [-k-1 + t_{0},t_{0} + 1]$ with
$|p-q|< \delta, |x(p)-x(q)|< \epsilon$. Take $\delta^{*}=\min(1,
\delta)$ and consider $s$ with $|t_{0}-s|< \delta^{*}$. Clearly,
for any $\theta \in [-k, 0], t_{0}+ \theta$ and $s+\theta$ both
belong to $[-k-1+t_{0},t_{0}+1]$ and hence
\begin{equation*}
\sup\{|x(t_{0}+\theta)-x(s+\theta)|\hbox{\rm :}\  \theta \in [-k,
0]\}< \epsilon.
\end{equation*}
That is, given $\epsilon > 0$, there is a $\delta^{*}>0$ such that
\begin{equation*}
\|u(s)-u(t_{0})\|_{k}< \epsilon
\end{equation*}
whenever $|s-t_{0}|< \delta^{*}$.

Next, we claim that for a given $\epsilon > 0$ and $t_{0} \in [0,
j\tau_{1}]$, there exists $\delta > 0$ such that
$p_{k}(u(t_{0})-u(t))< \epsilon$ whenever $t_{0} \in [0,
j\tau_{1}]$ satisfies $|t-t_{0}|<\delta$.

By Proposition 1.6, and the definition of $F$, for a given
$\varphi \in F, G_{k}\hbox{\rm :}\  [0, k\tau_{1}]\rightarrow
l^{1}(\mathbb{R})$ defined as
$[G_{k}(s)]^{(i)}=b_{n(k)+(i-1)}\varphi (s-\tau_{n(k)+(i-1)})$ is
a continuous function. Consider
\begin{align*}
\sum_{i=n(k)}^{K} b_{i}[u(t_{0})-u(t)](s-\tau_{i})&=
\sum_{i=n(k)}^{n(k+j)-1} b_{i}[u(t_{0})-u(t)](s-\tau_{i})\\[.5pc]
&\quad\,+
\sum_{i=n(k+j)}^{K} b_{i}[u(t_{0})-u(t)](s-\tau_{i})\\[.5pc]
&=\sum_{i=n(k)}^{n(k+j)-1}
b_{i}[x(t_{0}+s-\tau_{i})-x(t+s-\tau_{i})]\\[.5pc]
&=\sum_{i=n(k+j)}^{K} b_{i}[\phi(t_{0}+s-\tau_{i})-\phi
(t+s\!-\tau_{i})]
\end{align*}
and hence
\begin{align*}
\sum_{i=n(k)}^{\infty} b_{i}[u(t_{0})-u(t)](s-\tau_{i})
&=\sum_{i=n(k)}^{n(k+j)-1}
b_{i}[x(t_{0}+s-\tau_{i})-x(t+s-\tau_{i})]\\[.5pc]
&+\sum_{i=n(k+j)}^{\infty} b_{i}[\phi(t_{0}+
s-\tau_{i})-\phi(t+s\!-\tau_{i})].
\end{align*}
By the uniform continuity of $G_{n(k+j)}$ on $[0, (k+j)\tau_{1}]$,
given $\epsilon > 0$, there exists a $\delta_{1}>0$ such that
$|p-q|<\delta_{1}$ implies that
\begin{equation*}
\sum_{i=n(k+j)}^{\infty} |b_{i}| |[\phi
(p-\tau_{i})-\phi(q-\tau_{i})]|< \epsilon/2.
\end{equation*}
So, for $t_{0}$ and $t \in [0, j\tau_{1}]$ with $|t_{0}-t|<
\delta_{1}$,
\begin{equation*}
\sup\left\{\sum_{i=n(k+j)}^{\infty} |b_{i}|
|[\phi(p-\tau_{i})-\phi(q-\tau_{i})]|\hbox{\rm :}\  s\in
[0,k\tau_{1}]\right\}< \epsilon/2.
\end{equation*}
Since the expression $\sum_{i=n(k)}^{n(k+j)-1}
b_{i}[x(t_{0}+s-\tau_{i})-x(t+s-\tau_{i})]$ involves evaluation of
$x$ over a compact set, given $\epsilon>0$, there exists
$\delta_{2}>0$ such that
\begin{equation*}
\sum_{i=n(k)}^{n(k+j)-1}\sup\{|b_{i}|
|[x(t_{0}+s-\tau_{i})-x(t+s-\tau_{i})]|\}< \epsilon/2
\end{equation*}
whenever $t_{0}, t \in [0, j\tau_{1}]$ and $|t_{0}-t|<
\delta_{2}$. Our claim is now proved by taking
$\delta=\min(\delta_{1}, \delta_{2})$.

\setcounter{theore}{2}
\begin{pot}
{\rm By Lemma 2.4, for a given $\phi \in F$, we have a solution
$x_{\phi}\hbox{\rm :}\ \mathbb{R}\rightarrow \mathbb{R}$ to (1.1).

Define $S_{t}\hbox{\rm :}\  F\rightarrow F$ as
\begin{align*}
[S_{t}\phi](\theta) &=x_{\phi}(t+\theta),\quad t+\theta>0\\[.2pc]
&=\phi (t+\theta),\quad t+\theta\leq 0.
\end{align*}

Lemma 2.5 shows that $S_{t}\phi \in F$ for all $t\geq 0$ and the
map $u(t)=S_{t}\phi$ is a  continuous function from $[0, \infty)$
into $F$. From the definition of $S_{t}$ one can verify linearity
of each $S_{t}$ and the property $S_{t+s}=S_{t}S_{s}$. We need
only to check that each $S_{t}$ is a bounded linear map on $F$.
This follows from the estimates (2.9), 2.10) and (2.1)

Next, we have
\begin{align}
[u(t)](\theta)&=\phi (t+\theta),\ \  t+\theta \leq 0\nonumber\\[.4pc]
[u(t)](\theta)&=\phi (0)+ \int_{0}^{t + \theta} L(u(s))\hbox{d}s,\
\  t + \theta > 0 
\end{align}
and using (2.11), it is easy to see that $u$ is indeed a mild
solution to (1.3).

The Banach space ${\bf C}_{g}$ defined below was the phase space
used in [4] for the study of infinite delay equations of which
(1.1) is a special case. An interesting observation is that $F$
contains ${\bf C}_{g}$.}
\end{pot}

\setcounter{theore}{5}
\begin{exampl}{\rm
Let $g\hbox{\rm :}\  (-\infty, 0]\rightarrow [1, \infty)$ be a
continuous non-increasing function such that
\begin{equation*}
\sum_{i=1}^{\infty} b_{i} g(-\tau_{i}) < \infty. \tag{2.12}
\end{equation*}
Consider the space ${\bf C}_{g}$ of all $\varphi \in \hbox{\bf
C}(-\infty, 0]\rightarrow \mathbb{R}$ such that
\begin{equation*}
\|\varphi\|_{g} \mathop{=}\limits_{}^{{\rm
def}}\sup\left\{\frac{|\varphi(\theta)|}{g(\theta)}\hbox{\rm :}\
\theta \in (-\infty,0]\right\}< \infty.
\end{equation*}
Then ${\bf C}_{g}\subset F$.}
\end{exampl}

\begin{proof}
Consider
$|b_{i}\|\varphi(s-\tau_{i})|=|b_{i}|\frac{|\varphi(s-\tau_{i})|}{g(s-\tau_{i})}g(s-\tau_{i})$.
For $s \in [0, k\tau_{1}]$ and $i\geq n(k), 0\geq s-\tau_{i}\geq
-\tau_{i}$ and hence
\begin{equation*}
\sup\{|b_{i}\|\varphi(s-\tau_{i})|\hbox{\rm :}\  s \in [0,
k\tau_{1}]\}\leq (\|\phi\|_{g}) |b_{i}|g(-\tau_{i}).
\end{equation*}
Now, the assertion follows from (2.12).
\end{proof}

\section*{Acknowledgement}

Major part of this work was done while the author was with the
Department of Mathematics, Indian Institute of Science, Bangalore.
The author acknowledges the referee whose comments greatly
improved the presentation of the paper and who pointed out many
misprints.

\end{document}